




%

%
 \font\twelvebf=cmbx12
 \font\twelvett=cmtt12
 \font\twelveit=cmti12
 \font\twelvesl=cmsl12
 \font\twelverm=cmr12		\font\ninerm=cmr9
 \font\twelvei=cmmi12		\font\ninei=cmmi9
 \font\twelvesy=cmsy10 at 12pt	\font\ninesy=cmsy9
 \skewchar\twelvei='177		\skewchar\ninei='177
 \skewchar\seveni='177	 	\skewchar\fivei='177
 \skewchar\twelvesy='60		\skewchar\ninesy='60
 \skewchar\sevensy='60		\skewchar\fivesy='60
%
%

%
 \font\fourteenrm=cmr12 scaled 1200
 \font\seventeenrm=cmr12 scaled 1440
 \font\fourteenbf=cmbx12 scaled 1200
 \font\seventeenbf=cmbx12 scaled 1440
%
%

%
%
%
\font\tenmsb=msbm10
\font\twelvemsb=msbm10 scaled 1200
\newfam\msbfam

%
\font\tensc=cmcsc10
\font\twelvesc=cmcsc10 scaled 1200
\newfam\scfam

%
\def\seventeenpt{\def\rm{\fam0\seventeenrm}%
 \textfont\bffam=\seventeenbf	\def\bf{\fam\bffam\seventeenbf}}
\def\fourteenpt{\def\rm{\fam0\fourteenrm}%
 \textfont\bffam=\fourteenbf	\def\bf{\fam\bffam\fourteenbf}}
\def\twelvept{\def\rm{\fam0\twelverm}%
 \textfont0=\twelverm	\scriptfont0=\ninerm	\scriptscriptfont0=\sevenrm
 \textfont1=\twelvei	\scriptfont1=\ninei	\scriptscriptfont1=\seveni
 \textfont2=\twelvesy	\scriptfont2=\ninesy	\scriptscriptfont2=\sevensy
 \textfont3=\tenex	\scriptfont3=\tenex	\scriptscriptfont3=\tenex
 \textfont\itfam=\twelveit	\def\it{\fam\itfam\twelveit}%
 \textfont\slfam=\twelvesl	\def\sl{\fam\slfam\twelvesl}%
 \textfont\ttfam=\twelvett	\def\tt{\fam\ttfam\twelvett}%
 \scriptfont\bffam=\tenbf 	\scriptscriptfont\bffam=\sevenbf
 \textfont\bffam=\twelvebf	\def\bf{\fam\bffam\twelvebf}%
 \textfont\scfam=\twelvesc	\def\sc{\fam\scfam\twelvesc}%
 \textfont\msbfam=\twelvemsb	
 \baselineskip 14pt%
 \abovedisplayskip 7pt plus 3pt minus 1pt%
 \belowdisplayskip 7pt plus 3pt minus 1pt%
 \abovedisplayshortskip 0pt plus 3pt%
 \belowdisplayshortskip 4pt plus 3pt minus 1pt%
 \parskip 3pt plus 1.5pt
 \setbox\strutbox=\hbox{\vrule height 10pt depth 4pt width 0pt}}
\def\tenpt{\def\rm{\fam0\tenrm}%
 \textfont0=\tenrm	\scriptfont0=\sevenrm	\scriptscriptfont0=\fiverm
 \textfont1=\teni	\scriptfont1=\seveni	\scriptscriptfont1=\fivei
 \textfont2=\tensy	\scriptfont2=\sevensy	\scriptscriptfont2=\fivesy
 \textfont3=\tenex	\scriptfont3=\tenex	\scriptscriptfont3=\tenex
 \textfont\itfam=\tenit		\def\it{\fam\itfam\tenit}%
 \textfont\slfam=\tensl		\def\sl{\fam\slfam\tensl}%
 \textfont\ttfam=\tentt		\def\tt{\fam\ttfam\tentt}%
 \scriptfont\bffam=\sevenbf 	\scriptscriptfont\bffam=\fivebf
 \textfont\bffam=\tenbf		\def\bf{\fam\bffam\tenbf}%
 \textfont\scfam=\tensc		\def\sc{\fam\scfam\tensc}%
 \textfont\msbfam=\tenmsb	
 \baselineskip 12pt%
 \abovedisplayskip 6pt plus 3pt minus 1pt%
 \belowdisplayskip 6pt plus 3pt minus 1pt%
 \abovedisplayshortskip 0pt plus 3pt%
 \belowdisplayshortskip 4pt plus 3pt minus 1pt%
 \parskip 2pt plus 1pt
 \setbox\strutbox=\hbox{\vrule height 8.5pt depth 3.5pt width 0pt}}

%
\def\twelvepoint{%
 \def\small{\tenpt\rm}%
 \def\normal{\twelvept\rm}%
 \def\large{\fourteenpt\rm}%
 \def\huge{\seventeenpt\rm}%
 \footline{\hss\twelverm\folio\hss}%
 \normal}
\def\tenpoint{%
 \def\small{\tenpt\rm}%
 \def\normal{\tenpt\rm}%
 \def\large{\twelvept\rm}%
 \def\huge{\fourteenpt\rm}%
 \footline{\hss\tenrm\folio\hss}%
 \normal}

\tenpoint

%

%
\catcode`\@=11
%
%
\def\footnote#1{\edef\@sf{\spacefactor\the\spacefactor}#1\@sf
 \insert\footins\bgroup\small
 \interlinepenalty100	\let\par=\endgraf
 \leftskip=0pt		\rightskip=0pt
 \splittopskip=10pt plus 1pt minus 1pt	\floatingpenalty=20000
 \smallskip\item{#1}\bgroup\strut\aftergroup\@foot\let\next}
%
%
%
%
\def\hexnumber@#1{\ifcase#1 0\or 1\or 2\or 3\or 4\or 5\or 6\or 7\or 8\or
 9\or A\or B\or C\or D\or E\or F\fi}
\edef\msbfam@{\hexnumber@\msbfam}
\def\Bbb#1{\fam\msbfam\relax#1}
%
%
%
\catcode`\@=12

\twelvepoint

\font\twelvebf=cmbx12	
\font\ninebf=cmbx9	\font\sevenbf=cmbx7	\font\fivebf=cmbx5

\font\twelvebfit=cmbxti10 at 12pt	
\newfam\bfitfam

\font\twelvebm=cmmib10 at 12pt		\font\tenbm=cmmib10
\font\ninebm=cmmib9	\font\sevenbm=cmmib7	\font\fivebm=cmmib5

\skewchar\twelvebm='177	\skewchar\tenbm='177
\skewchar\ninebm='177	\skewchar\sevenbm='177	\skewchar\fivebm='177

\def\twelvepointbold{\def\bold{
\textfont0=\twelvebf	\scriptfont0=\ninebf	\scriptscriptfont0=\sevenbf
\textfont1=\twelvebm	\scriptfont1=\ninebm	\scriptscriptfont1=\sevenbm
\textfont\bffam=\twelvebf	\textfont\bfitfam=\twelvebfit
\def\rm{\fam\bffam\twelvebf}%
\def\it{\fam\bfitfam\twelvebfit}%
\rm}}

\twelvepointbold

\chardef\tempcat=\the\catcode`\@
\catcode`\@=11


\def\cydot{{\mathsurround=0pt$\cdot$}}

\def\ubar#1{\oalign{#1\crcr\hidewidth
    \vbox to.2ex{\hbox{\char22}\vss}\hidewidth}}

\def\cprime{\/{\mathsurround=0pt$'$}}
\def\Cprime{{\mathsurround=0pt$'$}}
\def\cdprime{\/{\mathsurround=0pt$''$}}
\def\Cdprime{{\mathsurround=0pt$\ubar{\hbox{$''$}}$}}

\def\dbar{dj}           
\def\Dbar{Dj}           

\def\dz{dz}
\def\Dz{Dz}
\def\dzh{dzh\cydot }
\def\Dzh{Dzh\cydot }


\def\@gobble#1{}
\def\@testgrave{\`}
\def\@stressit{\futurelet\chartest\@stresschar }

\def\@stresschar#1{%
  \ifx #1y\def\result{\futurelet\chartest\@yligature}%
  \else \ifx #1Y\def\result{\futurelet\chartest\@Yligature}%
  \else \ifx\chartest\@testgrave \def\result{\accent"26 }%
  \else \def\result{\accent"26 #1}%
  \fi \fi \fi
  \result }

\def\@yligature{%
  \ifx a\chartest \def\result{\accent"26 \char"1F \@gobble}%
  \else \ifx u\chartest \def\result{\accent"26 \char"18 \@gobble}%
  \else \def\result{\accent"26 y}%
  \fi \fi
  \result }

\def\@Yligature{%
  \ifx a\chartest \def\result{\accent"26 \char"17 \@gobble}%
  \else \ifx A\chartest \def\result{\accent"26 \char"17 \@gobble}%
  \else \ifx u\chartest \def\result{\accent"26 \char"10 \@gobble}%
  \else \ifx U\chartest \def\result{\accent"26 \char"10 \@gobble}%
  \else \def\result{\accent"26 Y}%
  \fi \fi \fi \fi
  \result }

\def\!{\ifmmode \mskip-\thinmuskip \fi}


\def\cyracc{%
  \def\cydot{{\kern0pt}}%
  \def\cprime{\char"7E }\def\Cprime{\char"5E }%
  \def\cdprime{\char"7F }\def\Cdprime{\char"5F }%
  \def\dbar{dj}\def\Dbar{Dj}%
  \def\dz{\char"1E }\def\Dz{\char"16 }%
  \def\dzh{\char"0A }\def\Dzh{\char"02 }%
  \def\'##1{\if c##1\char"0F %
    \else \if C##1\char"07 %
    \else \accent"26 ##1\fi \fi }%
  \def\`##1{\if e##1\char"0B %
    \else \if E##1\char"03 %
    \else \errmessage{accent \string\` not defined in cyrillic}%
        ##1\fi \fi }%
  \def\=##1{\if e##1\char"0D %
    \else \if E##1\char"05 %
    \else \if \i##1\char"0C %
    \else \if I##1\char"04 %
    \else \errmessage{accent \string\= not defined in cyrillic}%
        ##1\fi \fi \fi \fi }%
  \def\u##1{\if \i##1\accent"24 i%
    \else \accent"24 ##1\fi }%
  \def\"##1{\if \i##1\accent"20 \char"3D %
    \else \if I##1\accent"20 \char"04 %
    \else \accent"20 ##1\fi \fi }%
  \def\!{\ifmmode \def\result{\mskip-\thinmuskip}%
    \else \def\result{\@stressit}\fi \result}}


\def\keep@cyracc{\let\cyr=\relax \let\i=\relax
        \let\ubar=\relax \let\cydot=\relax
        \let\cprime=\relax \let\Cprime=\relax
        \let\cdprime=\relax \let\Cdprime=\relax
        \let\dbar=\relax \let\Dbar=\relax
        \let\dz=\relax \let\Dz=\relax
        \let\dzh=\relax \let\Dzh=\relax
        \let\'=\relax \let\`=\relax \let\==\relax
        \let\u=\relax \let\"=\relax \let\!=\relax }

\catcode`\@=\tempcat
\newfam\cyrfam
\font\twelvecyr=wncyr10 at 12pt
\def\cyr{\fam\cyrfam\twelvecyr\cyracc}

\newcount\EQNO      \EQNO=0
\newcount\FIGNO     \FIGNO=0
\newcount\REFNO     \REFNO=0
\newcount\SECNO     \SECNO=0
\newcount\SUBSECNO  \SUBSECNO=0
\newcount\FOOTNO    \FOOTNO=0
\newbox\FIGBOX      \setbox\FIGBOX=\vbox{}
\newbox\REFBOX      \newbox\REFBOXTMP
\setbox\REFBOX=\vbox{\bigskip\centerline{\bf REFERENCES}\smallskip}
\newbox\Partialpage \newdimen\Mark
\newdimen\REFSIZE   \REFSIZE=\vsize

\expandafter\ifx\csname normal\endcsname\relax\def\normal{\null}\fi

\def\MultiRef#1{\global\advance\REFNO by 1 \nobreak\the\REFNO%
    \global\setbox\REFBOX=\vbox{\unvcopy\REFBOX\normal
      \smallskip\item{\the\REFNO .~}#1}%
    \gdef\label##1{\xdef##1{\nobreak[\the\REFNO]}}%
    \gdef\Label##1{\xdef##1{\nobreak\the\REFNO}}}
\def\NoRef#1{\global\advance\REFNO by 1%
    \global\setbox\REFBOX=\vbox{\unvcopy\REFBOX\normal
      \smallskip\item{\the\REFNO .~}#1}%
    \gdef\label##1{\xdef##1{\nobreak[\the\REFNO]}}}
\def\Eqno{\global\advance\EQNO by 1 \eqno(\the\EQNO)%
    \gdef\label##1{\xdef##1{\nobreak(\the\EQNO)}}}
\def\Eqalignno{\global\advance\EQNO by 1 &(\the\EQNO)%
    \gdef\label##1{\xdef##1{\nobreak(\the\EQNO)}}}
\def\Fig#1{\global\advance\FIGNO by 1 Figure~\the\FIGNO%
    \global\setbox\FIGBOX=\vbox{\unvcopy\FIGBOX
      \narrower\smallskip\item{\bf Figure \the\FIGNO~~}#1}}
\def\Ref#1{\global\advance\REFNO by 1 \nobreak[\the\REFNO]%
    \global\setbox\REFBOX=\vbox{\unvcopy\REFBOX\normal
      \smallskip\item{\the\REFNO .~}#1}%
    \gdef\label##1{\xdef##1{\nobreak[\the\REFNO]}}}
\def\Section#1{\SUBSECNO=0\advance\SECNO by 1
    \bigskip\leftline{\bf \the\SECNO .\ #1}\nobreak}
\def\Subsection#1{\advance\SUBSECNO by 1
    \medskip\leftline{\bf \ifcase\SUBSECNO\or
    a\or b\or c\or d\or e\or f\or g\or h\or i\or j\or k\or l\or m\or n\fi
    )\ #1}\nobreak}
\def\Footnote#1{\global\advance\FOOTNO by 1 
    \footnote{\nobreak$\>\!{}^{\the\FOOTNO}\>\!$}{#1}
    \gdef\label##1{\xdef##1{$\>\!{}^{\the\FOOTNO}\>\!$}}}
\def\SameFootnote{$\>\!{}^{\the\FOOTNO}\>\!$}

\def\References{
     {\output={\global\setbox\Partialpage=\vbox{\unvbox255}}\eject}
     \Mark=\vsize
     \ifdim \ht\Partialpage > 0pt
          \advance\Mark by -\ht\Partialpage
          \advance\Mark by -0.1in
          \vbox{\unvbox\Partialpage}
     \fi
     \setbox\REFBOX=\vbox{\unvcopy\REFBOX\vfill\eject}
     \ifdim \ht\REFBOX > \Mark
          \setbox\REFBOXTMP=\vsplit\REFBOX to \Mark
          \vfill\copy\REFBOXTMP
     \else
          \REFSIZE=\Mark
     \fi
     \loop
          \setbox\REFBOXTMP=\vsplit\REFBOX to \REFSIZE
          \copy\REFBOXTMP
          \ifdim \ht\REFBOX > 0pt
     \repeat}
\def\NewRefPage{\setbox\REFBOX=\vbox{\unvcopy\REFBOX\vfill\eject}}


\newcount\LEMMA	  \LEMMA=0
\newcount\THM     \THM=0
\def\Lemma#1{\global\advance\LEMMA by 1\smallskip
    {\narrower\narrower\narrower\item{\bf Lemma~\the\LEMMA:}
    {\it #1}\smallskip}}
\def\Theorem#1{\global\advance\THM by 1\smallskip
    {\narrower\narrower\narrower\item{\bf Theorem~\the\THM:}
    {\it #1}\smallskip}}
\def\Heading#1#2{\smallskip\goodbreak
    {\narrower\narrower\narrower\item{\bf #1:}
    {\it #2}\smallskip}}

\input epsf
\def\Fig#1#2#3#4#5{\global\advance\FIGNO by 1 Figure~\the\FIGNO#5
    \topinsert
    \centerline{\epsfysize=#4\epsffile[#3]{#2}}
    {\bigskip\hsize=5.5in\hskip\parindent
     \vbox{\small\item{\bf Figure~\the\FIGNO:}{#1}}}
    \bigskip\endinsert}

\def\today{\number\day\space\ifcase\month\or
  January\or February\or March\or April\or May\or June\or
  July\or August\or September\or October\or November\or December\fi
  \space\number\year}

\def\RR{{\Bbb R}}
\def\CC{{\Bbb C}}
\def\HH{{\Bbb H}}
\def\OO{{\Bbb O}}
\def\AA{{\cal A}}

\def\AAA{{\Bbb A}}
\def\VVV{{\Bbb V}}

\def\bar{\overline}
\let\tilde=\widetilde

\def\Tr{{\rm tr\,}}
\def\Re{{\rm Re}}
\def\Im{{\rm Im}}

\def\eplus{\Psi^{\scriptscriptstyle+}}
\def\eminus{\Psi^{\scriptscriptstyle-}}

\def\Rule{3-$\Psi$'s rule}


\rightline{3 May 2000; revised 11 October 2000}
\bigskip

\centerline{\large\bf OCTONIONIC HERMITIAN MATRICES}
\smallskip
\centerline{\large\bf WITH NON-REAL EIGENVALUES}
\bigskip

\centerline{Tevian Dray}
\centerline{\it Department of Mathematics, Oregon State University,
		Corvallis, OR  97331, USA}
\centerline{\tt tevian{\rm @}math.orst.edu}
\medskip
\centerline{Jason Janesky}
\centerline{\it Department of Physics, Oregon State University,
		Corvallis, OR  97331, USA
\Footnote{Present address:
	Phoenix Corporate Research Laboratories,
	Motorola Inc.,
	Tempe, AZ  85284, USA,
	{\tt r47569{\rm @}email.sps.mot.com}
}}

\medskip
\centerline{Corinne A. Manogue}
\centerline{\it Department of Physics, Oregon State University,
		Corvallis, OR  97331, USA}
\centerline{\tt corinne{\rm @}physics.orst.edu}

\bigskip\bigskip
\centerline{\bf ABSTRACT}
\midinsert
\narrower\narrower\noindent
We extend previous work on the eigenvalue problem for Hermitian octonionic
matrices by discussing the case where the eigenvalues are not real, giving a
complete treatment of the $2\times2$ case, and summarizing some prelimenary
results for the $3\times3$ case.
\endinsert
\bigskip

\Section{INTRODUCTION}

In previous work
[\MultiRef{Tevian Dray and Corinne A. Manogue,
{\it The Octonionic Eigenvalue Problem},
Adv.\ Appl.\ Clifford Algebras {\bf 8}, 341--364 (1998).}%
\label\Eigen\Label\EIGEN
,\MultiRef{Tevian Dray and Corinne A. Manogue,
{\it Finding Octonionic Eigenvectors Using {\sl Mathematica}},
Comput.\ Phys.\ Comm.\ {\bf 115}, 536--547 (1998).}\label\Find\Label\FIND
; see also
\MultiRef{Susumu Okubo,
{\it Eigenvalue Problem for Symmetric $3\times3$ Octonionic Matrix},
Adv.\ Appl.\ Clifford Algebras {\bf 9}, 131--176 (1999).}]\label\Okubo
, we considered the real eigenvalue problem for $2\times2$ and $3\times3$
Hermitian matrices over the octonions $\OO$.  The $2\times2$ case corresponds
closely to the standard, complex eigenvalue problem, since any $2\times2$
octonionic Hermitian matrix lies in a complex subalgebra $\CC\subset\OO$.  The
$3\times3$ case requires considerable care, resulting in some changes in the
expected results.  However, we also showed in \Eigen\ that there are
octonionic Hermitian matrices which admit eigenvalues which are not real, and
which it is the purpose of this paper to discuss.
\Footnote{The quite different {\it Jordan\/} eigenvalue problem case admits
only real eigenvalues and was discussed in~%
\Ref{Tevian Dray and Corinne A. Manogue,
{\it The Exceptional Jordan Eigenvalue Problem},
Internat.\ J.\ Theoret.\ Phys.\ {\bf 38}, 2901--2916 (1999).}\label\Jordan
.}

We consider both the right eigenvalue problem
$$A v = v\lambda \Eqno$$\label\Master
and the left eigenvalue problem
$$A w = \lambda w \Eqno$$\label\Orig
where $A$ is a Hermitian octonionic matrix.  The $2\times2$ case is reasonably
straightforward, and can be completely solved.  Although we argue that the
right eigenvalue problem \Master\ is more fundamental, we obtain some
intriguing results relating the sets of left and right eigenvectors, as well
as the matrices whose eigenvectors they are.

\goodbreak
We then briefly discuss our preliminary results in the $3\times3$ case.
Although we have been able to obtain a 3rd-order characteristic equation for
the (right) eigenvalues in this case, we have not been able to solve this
equation, nor have we been able to extend our orthonormality results
[\EIGEN,\FIND] from the real case.

Both of these cases have applications to physics.  Three of the four
superstring equations of motion can be written as (2 separate) $2\times2$
octonionic eigenvalue problems with real eigenvalue $0$
[\MultiRef{David B.~Fairlie \& Corinne A.~Manogue,
{\it Lorentz Invariance and the Composite String},
Phys.\ Rev.\ {\bf D34}, 1832--1834 (1986)},%
\MultiRef{David B.~Fairlie \& Corinne A.~Manogue,
{\it A Parameterization of the Covariant Superstring},
Phys.\ Rev.\ {\bf D36}, 475-479 (1987).},\label\Fairlie
\MultiRef{Corinne A.~Manogue \& Anthony Sudbery,
{\it General Solutions of Covariant Superstring Equations of Motion},
Phys.\ Rev.\ {\bf D40}, 4073-4077 (1989).}]\label\Sudbery
.  The resulting equation is really the massless Dirac (Weyl) equation in 10
dimensions (and in momentum space), and this has recently been used in a model
for dimensional reduction~%
[\MultiRef{Corinne A. Manogue and Tevian Dray,
{\it Dimensional Reduction},
Mod.\ Phys.\ Lett.\ {\bf A14}, 93--97 (1999).}\label\Dim\Label\DIM
,\MultiRef{Corinne A. Manogue and Tevian Dray,
{\it Quaternionic Spin},
in {\bf Clifford Algebras and their Applications in Mathematical Physics},
eds.\ Rafa\l \ Ab\l amowicz and Bertfried Fauser,
Birkh\"auser, Boston, 2000, pp.\ 29--46.
\hfill{\small({\tt http://xxx.lanl.gov/abs/hep-th/9910010})}}]%
\label\Spin\Label\SPIN
.  In this model each quaternionic subalgebra \hbox{$\HH\subset\OO$} gives
rise to a spectrum of (free) particles which corresponds exactly with the
spins and helicities of a generation of leptons, including both massless and
massive particles.  Furthermore, in a natural sense there are precisely 3 such
quaternionic subalgebras compatible with the dimensional reduction, which we
have interpreted as generations.

Fundamental to this model is the use of the Lorentz group
$SL(2,\CC)\subset SL(2,\OO)$ to analyze the spin states of the resulting
particles
[\MultiRef{Corinne A. Manogue and J\"org Schray,
{\it Finite Lorentz transformations, automorphisms, and division algebras},
J. Math.\ Phys.\ {\bf 34}, 3746--3767 (1993).},%
\MultiRef{J\"org Schray,
{\bf Octonions and Supersymmetry},
Ph.D.\ thesis, Department of Physics, Oregon State University, 1994.}]%
\Label\Thesis
.  This is again an eigenvalue problem, this time for an octonionic
self-adjoint operator.  We discuss this operator eigenvalue problem below,
showing that already at the quaternionic level it admits eigenvalues which are
not real.  This leads to spin states which {\it are\/} simultaneous
eigenstates of all the spin operators, although not all the eigenvalues are
real.  This result could have implications for the interpretation of quantum
mechanics \Spin.

The $3\times3$ case is of particular interest mathematically because it
corresponds to the exceptional Jordan algebra, also known as the Albert
algebra.  There have been numerous attempts to use this algebra to describe
quantum physics, which was in fact Jordan's original motivation.  More
recently, Schray 
[\Thesis,%
\MultiRef{J\"org Schray,
{\it The General Classical Solution of the Superparticle},
Class.\ Quant.\ Grav.\ {\bf 13}, 27 (1996)}] \label\Schray\Label\SCHRAY
has shown how to use the exceptional Jordan algebra to give
an elegant description of the superparticle, which we have been attempting to
extend to the superstring.  Our dimensional reduction scheme extends naturally
to this case \Jordan, and we believe it is the natural language to describe
the fundamental particles of nature.

The paper is organized as follows.  In Section 2 we briefly review the
properties of octonions.  In Section 3 we consider $2\times2$ octonionic
Hermitian matrices, and in Section 4 we discuss the $2\times2$ self-adjoint
spin operators.  In Section 5, we summarize our preliminary attempts to
generalize these results to $3\times3$ octonionic Hermitian matrices.  Along
the way, we have need of several identities involving octonionic associators,
which are closely related to the ``\Rule'' needed for supersymmetric theories;
this is discussed in the Appendix.  Finally, in Section 6 we discuss our
results.

\Section{OCTONIONS}

We summarize here only the essential properties of the octonions $\OO$.  For a
more detailed introduction, see \Eigen\ or
[\MultiRef{Feza G\"ursey and Chia-Hsiung Tze,
{\bf On the Role of Division, Jordan, and Related Algebras in Particle Physics},
World Scientific, Singapore, 1996.}%
,\MultiRef{S. Okubo,
{\bf Introduction to Octonion and Other Non-Associative Algebras in Physics},
Cambridge University Press, Cambridge, 1995.}]%
.%

The octonions $\OO$ are the nonassociative, noncommutative, normed division
algebra over the reals.  In terms of a natural basis, an octonion $a$ can be
written
$$a = \sum\limits_{q=1}^8 a^q e_q \Eqno$$
where the coefficients $a^q$ are real, and where the basis vectors satisfy
$e_1=1$ and
$$e_q^2 = -1 \qquad (q=2,...,8) \Eqno$$
The multiplication table is
conveniently encoded in the 7-point projective plane, shown in
\Fig{The representation of the octonionic multiplication table using the
7-point projective plane, where we have used the conventional names
$\{i,j,k,k\ell,j\ell,i\ell,\ell\}$ for $\{e_2,...,e_8\}$.  Each of the 7
oriented lines gives a quaternionic triple.}
{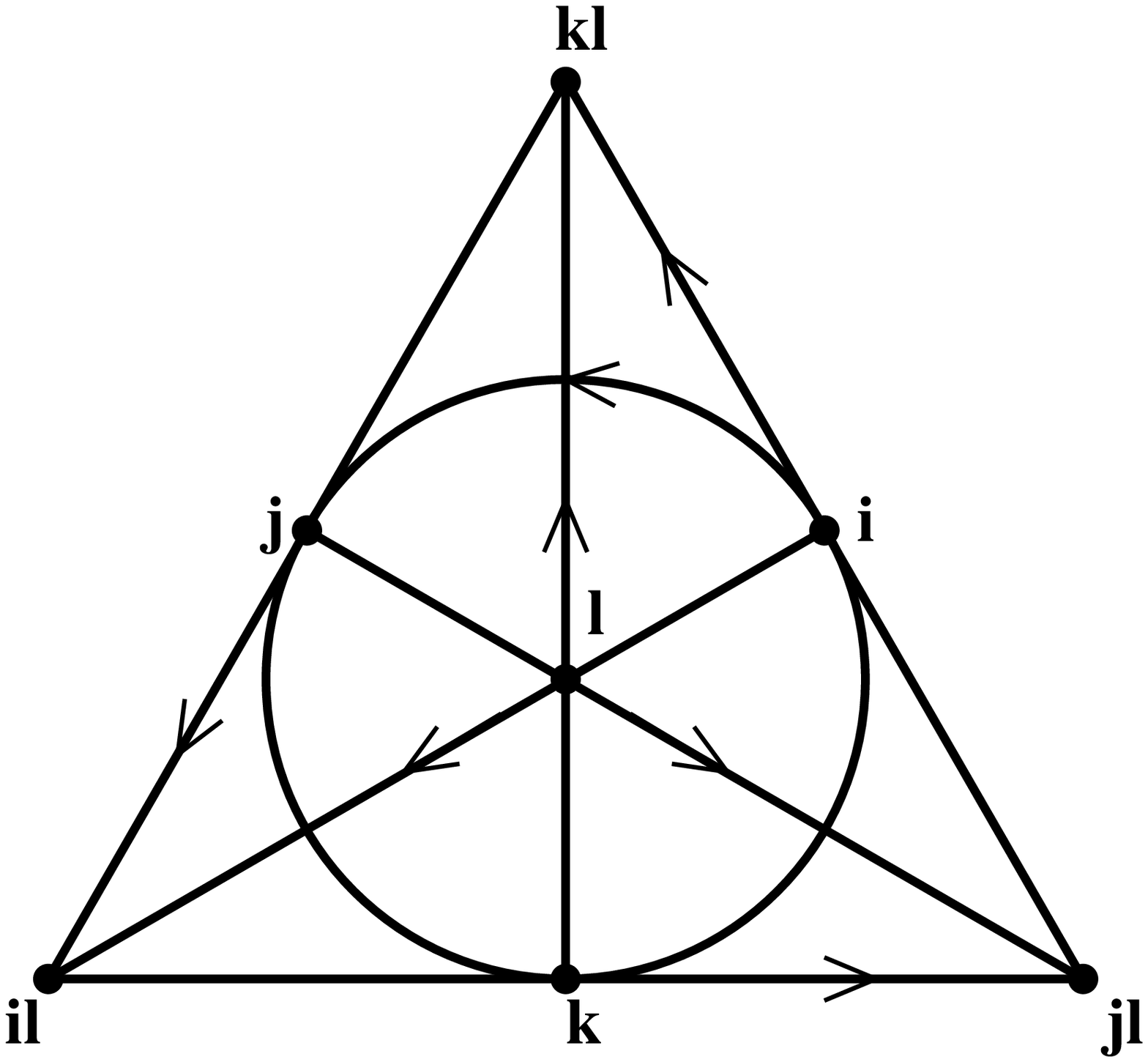}{68 168 543 614}{3in}
{.  The product of any two imaginary units is given by the third unit on the
unique line connecting them, with the sign determined by the relative
orientation.}

{\it Octonionic conjugation\/} is given by reversing the sign of the
imaginary basis units
$$\bar a = a^1 e_1 - \sum\limits_{q=2}^8 a^q e_q \Eqno$$
Conjugation is an antiautomorphism, since it satisfies
$$\bar{ab} = \bar{b} \> \bar{a}$$
The real and imaginary parts of an octonion $a$ are given by
$$\Re(a) = {1\over2} (a+\bar{a})
  \qquad\qquad
  \Im(a) = {1\over2} (a-\bar{a})
  \Eqno$$

The {\it inner product\/} on $\OO$ is the one inherited from $\RR^8$, namely
$$a \cdot b = \sum_q a^q b^q \Eqno$$
which can be rewritten as
$$a \cdot b
  = {1\over2} (a \bar{b} + b \bar{a})
  = {1\over2} (\bar{b} a + \bar{a} b)
  \Eqno$$
and which satisfies the identities
$$\eqalignno{
a \cdot (xb) &= b \cdot (\bar{x}a) \Eqalignno\label\DotID \cr
(ax) \cdot (bx) &= |x|^2 \> a \cdot b \Eqalignno\label\DotIDii
  }$$
for any $a,b,x\in\OO$.  The {\it norm\/} of an octonion is just
$$|a| = \sqrt{a \bar{a}} = \sqrt{a \cdot a} \Eqno$$
which satisfies the defining property of a normed division algebra, namely
$$|ab| = |a| |b| \Eqno$$

The {\it associator\/} of three octonions is
$$[a,b,c] = (ab)c - a(bc) \Eqno$$
which is totally antisymmetric in its arguments, has no real part, and changes
sign if any one of its arguments is replaced by its octonionic conjugate.
Although the associator does not vanish in general, the octonions do satisfy a
weak form of associativity known as {\it alternativity}, namely
$$[b,a,a]=0=[b,a,\bar{a}] \Eqno$$
The underlying reason for alternativity is Artin's Theorem
[\MultiRef{
Richard D. Schafer,
{\bf An Introduction to Nonassociative Algebras},
Academic Press, New York, 1966; reprinted by Dover, Mineola NY, 1995.},%
\MultiRef{Emil Artin,
{\bf Geometric Algebra},
Interscience Publishers, New York, 1957; reprinted by
John Wiley \& Sons, New York, 1988.}\label\Artin
]%
, which states that that any two octonions lie in a quaternionic subalgebra of
$\OO$, so that any product containing only two octonionic directions is
associative.  We will also have use for the associator identity
$$[a,b,c]d + a[b,c,d] = [ab,c,d] - [a,bc,d] + [a,b,cd] \Eqno$$\label\Assoc
for any $a,b,c,d\in\OO$, which is proved by writing out all the terms.

\Section{$\bold 2\times2$ OCTONIONIC HERMITIAN MATRICES}

The general $2\times2$ octonionic Hermitian matrix can be written
$$A = \pmatrix{p& a\cr \noalign{\smallskip} \bar{a}& m\cr} \Eqno$$
with $p,m\in\RR$ and $a\in\OO$, and satisfies its characteristic equation
$$A^2 - (\Tr A) \, A + (\det A) \, I = 0 \Eqno$$
where $\Tr A$ denotes the trace of $A$, and where there is no difficulty with
commutativity and associativity in defining the determinant of $A$ as usual
via
$$\det A = pm - |a|^2 \Eqno$$
since the components of $A$ lie in a complex subalgebra $\CC\subset\OO$.  If
$a=0$ the eigenvalue problem is trivial, so we assume $a\ne0$.  We also set
$$v=\pmatrix{x\cr y\cr} \Eqno$$\label\vform
with $x,y\in\OO$.
\goodbreak

\Subsection{Left Eigenvalue Problem}

As pointed out in \Eigen, even quaternionic Hermitian matrices can admit left
eigenvalues which are not real, as is shown by the following example:
$$\pmatrix{1&-i\cr \noalign{\smallskip} i&~~1}
	\pmatrix{1\cr \noalign{\smallskip} k}
  = \pmatrix{1+j\cr \noalign{\smallskip} k+i}
  = (1+j) \pmatrix{1\cr \noalign{\smallskip} k}
  \Eqno$$\label\Unex

Direct computation allows us to determine which Hermitian matrices $A$ admit
left eigenvalues which are not real.  Inserting \vform\ in \Orig\ leads to
$$(\lambda-p) x = a y \qquad\qquad (\lambda-m) y = \bar{a} x \Eqno$$
which in turn leads to
$${\bar{a} \Bigl((\lambda-p) x\Bigr) \over |a|^2}
  = y = {(\bar\lambda-m) (\bar{a} x) \over |\lambda-m|^2}
  \Eqno$$
Assuming without loss of generality that $x\ne0$ and taking the norm of both
sides yields
$$|a|^2 = |\lambda-p| |\lambda-m| \Eqno$$\label\Norm
resulting in
$${\bar{a} \Bigl((\lambda-p) x\Bigr) \over |\lambda-p|}
   = {(\bar\lambda-m) (\bar{a} x) \over |\lambda-m|}
  \Eqno$$\label\EqLeft
This equation splits into two independent parts, the terms (in the numerator)
which involve the imaginary part of $\lambda$, which is nonzero by assumption,
\Footnote{We can in fact assume without loss of generality that
$\Re(\lambda)=0$ by replacing $A$ with $A-\Re(\lambda)I$.}\label\wlog
and those which don't.  Looking first at the latter leads to
$$p = m \Eqno$$
which in turn reduces \EqLeft\ to
$$\bar{a} (\lambda x) = \bar\lambda (\bar{a} x) \Eqno$$
which forces $a$ to be purely imaginary (and $\lambda\cdot a=0$), but
which puts no conditions on $x$.

Denoting the $2\times2$ identity matrix by  $I$ and setting
$$J(\hat{r})
  = \pmatrix{0& -\hat{r}\cr \noalign{\smallskip} \hat{r}& ~~0\cr} \Eqno$$
for any pure imaginary unit octonion $\hat{r}$, and noting that this latter
condition can be written as $\hat{r}^2=-1$, we have
\Lemma{The set of $2\times2$ Hermitian matrices $A$ for which left
eigenvalues exist which are not real is
$$\AAA := \{ A{:}~ A = p \, I + q \, J(\hat{r}) ; \quad
  p,q\in\RR, \, q\ne0, \, \hat{r}^2=-1 \} \Eqno$$\label\Adef
}

The set $\AAA$ has some remarkable properties, which will be further discussed
below.  Without loss of generality, we can take $\hat{r}=i$, so that $A$ takes
the form
$$A = \pmatrix{p&-iq\cr \noalign{\smallskip} iq&~~p} \Eqno$$\label\Aform
Let us find the general solution of the left eigenvalue problem for these
matrices.  Taking $A$ as in \Aform\ and $v$ as in \vform\ we can rewrite
\Orig\ as
$${\lambda - p \over q} \> x = -iy \qquad
  {\lambda - p \over q} \> y = ix 
  \Eqno$$\label\Leq
Taking the norm of both sides immediately yields
$$|x|^2 = |y|^2 \Eqno$$\label\xyNorm
and we can normalize both of these to 1 without loss of generality.  We thus
obtain
$$\eqalign{{\lambda - p \over q}
  &= -(iy)\bar{x} = (ix)\bar{y} \cr
  &= -[i,y,\bar{x}] - i(y\bar{x}) = [i,x,\bar{y}] + i(x\bar{y}) \cr
  }\Eqno$$
But since
$$[z,y,\bar{x}] = -[z,y,x] = [z,x,y] = -[z,x,\bar{y}] \Eqno$$
for any $z$, the two associators cancel, and we are left with
$$x \cdot y = 0 \Eqno$$\label\xyDot
Thus, $x$ and $y$ correspond to orthonormal vectors in $\OO$ thought of as
$\RR^8$.  This argument is fully reversible; any suitably normalized $x$ and
$y$ which are orthogonal yield an eigenvector of $A$.  We have therefore shown
that {\it all\/} matrices in $\AAA$ have the same left eigenvectors:

\Lemma{The set of left eigenvectors for any matrix $A\in\AAA$ is given by
$$\VVV = \left\{
  \pmatrix{x\cr y}{:}~ |x|^2=|y|^2 ;~ x \cdot y =0 \right\} \Eqno$$
}

\noindent
The left eigenvalue is given in each case by \Leq.  Furthermore, left
multiplication by an arbitrary octonion preserves the set $\VVV$, so that
matrices in $\AAA$ have the property that left multiplication of left
eigenvectors yields another left eigenvector (albeit with a different
eigenvalue).
\Footnote{Direct computation shows that, other than real matrices, the
matrices in $\AAA$ are the only $2\times2$ Hermitian matrices with this
property.}
It follows from \Leq\ and \xyNorm\ that
$$|\lambda - p| = q \Eqno$$\label\pqeq
Inserting this into either of \Leq, multiplying both sides by $i$, and using
the identities \DotID\ and \DotIDii\ then shows that \xyDot\ forces
$$\lambda\cdot i = 0\Eqno$$
However, these are the only restrictions on $\lambda$, since \Leq\ can be used
to construct eigenvectors having {\it any\/} eigenvalue satisfying these two
conditions.

\Subsection{Right Eigenvalue Problem}

As discussed in \Eigen, the {\it right} eigenvalues of quaternionic Hermitian
matrices must be real, which is a strong argument in favor of \Master\ over
\Orig.  However, as pointed out in \Eigen, there do exist {\it  octonionic}
Hermitian matrices which admit right eigenvalues which are not real, as is
shown by the following example:
$$\pmatrix{1&-i\cr \noalign{\smallskip} i&~~1}
	\pmatrix{j\cr \noalign{\smallskip} \ell}
  = \pmatrix{j-i\ell\cr \noalign{\smallskip} \ell+k}
  = \pmatrix{j\cr \noalign{\smallskip} \ell} (1+k\ell)
  \Eqno$$\label\SpecEx

Proceeding as we did for left eigenvectors, we can determine which matrices
$A$ admit right eigenvalues which are not real.  Inserting \vform\ into
\Master\ leads to
$$x (\lambda-p) = a y \qquad\qquad y (\lambda-m) = \bar{a} x \Eqno$$
\label\Right
which in turn leads to
$${\bar{a} \Bigl(x (\lambda-p)\Bigr) \over |a|^2}
  = y = {(\bar{a} x) (\bar\lambda-m) \over |\lambda-m|^2}
  \Eqno$$
Taking the norm of both sides (and assuming $x\ne0$) again yields \Norm,
resulting in
$${\bar{a} \Bigl(x (\lambda-p)\Bigr) \over |\lambda-p|}
   = {(\bar{a} x) (\bar\lambda-m) \over |\lambda-m|}
  \Eqno$$\label\EqRight
Just as for the left eigenvector problem, this equation splits into two
independent parts, the terms (in the numerator) which involve the imaginary
part of $\lambda$, which is nonzero by assumption, and those which don't.\wlog\
Looking first at the latter again forces $p=m$, which in turn forces
$|y|=|x|$.  The remaining condition is now
$$\bar{a} (x \lambda) = (\bar{a} x) \bar\lambda \Eqno$$
so that $\bar{a}$, $\Im(\lambda)$, and $x$ antiassociate.  In particular, this
forces both $a$ and $x$ to be pure imaginary, as well as
$$\eqalignno{
\lambda \cdot a &= 0 \Eqalignno\label\laz \cr
\lambda\cdot x &= 0 = a\cdot x \Eqalignno\label\Xdot
  }$$
with corresponding identities also holding for $y$.
\Footnote{This also implies that $a\lambda\cdot x=0$, that is, $x$ (and $y$)
must be orthogonal to the quaternionic subalgebra generated by $\lambda$ and
$a$.}
We conclude that the necessary and sufficient condition for matrices to admit
right eigenvalues which are not real is that $A\in\AAA$:

\Lemma{The set of $2\times2$ Hermitian matrices $A$ for which right
eigenvalues exist which are not real is $\AAA$ as defined in \Adef.}

Thus, all $2\times2$ Hermitian matrices which admit right eigenvalues which
are not real also admit left eigenvalues which are not real, and vice versa!

\Heading{Corollary}
{A $2\times2$ octonionic Hermitian matrix admits right eigenvalues which are
not real if and only if it admits left eigenvalues which are not real.}

Turning to the eigenvectors, inserting $p=m$ into \Right\ leads to
$$\bar{x} (ay) = \bar{y} (\bar{a}x) \Eqno$$
and inserting the conditions on $a$, $x$, and $y$ now leads to
$$x \cdot y = 0 \Eqno$$
just as for left eigenvectors.  All right eigenvectors with non-real
eigenvalues are hence in $\VVV$, although the converse is false (since right
eigenvectors have no real part).  Furthermore, not all of the remaining
elements of $\VVV$ will be eigenvectors for any given matrix $A$ (since right
eigenvectors have no ``quaternionic'' part).

Putting all of this together, typical solutions of the (right) eigenvalue
problem for $A$ as in \Aform\ can thus be written as
$$\eqalign{
  v &= n \pmatrix{j\cr k\bar{s}\cr} \qquad
    \lambda_v = p+q\bar{s} \cr
\noalign{\smallskip}
  w &= n \pmatrix{k\bar{s}\cr j\cr} \qquad
    \lambda_w = p-q\bar{s} \cr
  }\Eqno$$\label\GenEx
where $p,q,n\in\RR$ and where
$$s=\cos\theta+k\ell\sin\theta \Eqno$$\label\Seq
The example given in \SpecEx\ is a special case of the first of \GenEx\ with
$p=q=n=1$ and $\theta=\pi/2$.

\Subsection{Characteristic Equation}

We now derive a generalized characteristic equation which is satisfied by
right eigenvalues of $2\times2$ octonionic Hermitian matrices.
\Footnote{Analogous results hold for left eigenvalues, but they are much less
elegant.}
Along the way, we also rederive some of the results of the previous
subsection.

Multiplying the first of \Right\ on the left by $\bar{a}$ and the second of
\Right\ on the right by $\lambda$ and subtracting leads to
$$y \left( \lambda^2 - \lambda (\Tr A) + (\det A) \right)
  = [a,y,\lambda] \Eqno$$
which can be solved for the characteristic equation in the form
$$\lambda^2 - \lambda (\Tr A) + (\det A)
  = {\bar{y} \, [\bar{a},x,\lambda] \over |y|^2} \Eqno$$ \label\CharI
Using \Assoc, we have
$$\bar{y} \, [\bar{a},x,\lambda]
  = [\bar{y}\,\bar{a},x,\lambda] - [\bar{y},\bar{a}x,\lambda]
     + [\bar{y},\bar{a},x\lambda] - [\bar{y},\bar{a},x] \, \lambda
  \Eqno$$\label\yaxl
But \Right\ immediately implies
$$[ay,x,\lambda] = 0 = [y,\bar{a}x,\lambda] \Eqno$$
so that the first 2 terms on the right-hand-side of \yaxl\ vanish.  Using
\Right\ again brings \yaxl\ to the form
$$\bar{y} \, [\bar{a},x,\lambda]
  = [\bar{y},\bar{a},xp+ay] - [\bar{y},\bar{a},x]\lambda
  = [\bar{y},\bar{a},x] \, (p-\lambda)
  \Eqno$$
and inserting this into \CharI\ leads finally to the {\it generalized
characteristic equation} for $\lambda$, namely
$$\lambda^2 - \lambda (\Tr A) + (\det A)
  = [\bar{a},x,y] \, {(\lambda-p) \over |y|^2}
  = [a,y,x] \, {(\lambda-m) \over |x|^2}
  \Eqno$$\label\CharII
where the final equality follows by symmetry.

If the associator $[a,x,y]$ vanishes, then $\lambda$ satisfies the ordinary
characteristic equation, and hence is real (since $A$ is complex Hermitian).
Otherwise, comparing real and imaginary parts of the last two terms in
\CharII\ provides an alternate derivation of $|y| = |x|$, and we recover $p=m$
as expected.  Furthermore, since the left-hand-side of \CharII\ lies in a
complex subalgebra of $\OO$, so does the right-hand-side, and it is then
straightforward to solve for $\lambda$ by considering its real and imaginary
parts.  The generalized characteristic equation \CharII\ then yields the
following equation for $\lambda$
$$\Bigl(\Re(\lambda)\Bigr)^2 - \Re(\lambda) (\Tr A) + (\det A)
  = \Bigl(\Im(\lambda)\Bigr)^2 < 0
  \Eqno$$\label\CharIIa
together with the requirement that
$${[\bar{a},x,y] \over |x| |y|} = 2 \, \Im(\lambda) \Eqno$$\label\ALT
The explicit form of the eigenvalues given in \GenEx\ and \Seq\ verifies that
there are no further restrictions on $\lambda$ other than \laz\ and \CharIIa.
Furthermore, having shown in the previous subsection that $a$ and $x$ (and
therefore also $y$) are pure imaginary, \ALT\ yields an alternate derivation
that $\lambda$ is orthogonal to $a$, which is \laz, as well as to $x$ and $y$,
which is \Xdot.

It follows directly from the generalized characteristic equation \CharIIa\
that
$$|\Im(\lambda)|\le|a| \Eqno$$
but this can be made more precise.  Given only that $m=p$, it is
straightforward to rewrite \CharIIa\ as (compare \pqeq)
$$|\lambda-p|^2 = |a|^2 \Eqno$$\label\DetRule
and it is intriguing that this seems to be almost the condition for the
vanishing of the determinant of $Q=A-\lambda I$.  But if $\lambda\not\in\RR$,
$Q$ is not Hermitian, and hence has no well-defined determinant.  However,
$QQ^\dagger$ is Hermitian, and $\det(QQ^\dagger)=0$ does indeed reduce to
\DetRule\ for $A$ of the form \Aform, that is with $p=m$ and $\Re(a)=0$,
provided that in addition \laz\ is assumed to hold.
\Footnote{Using the (square root of the) determinant of the Hermitian square
$QQ^\dagger$ is just Dieudonn\'e's prescription \Artin\ for the determinant of
a $2\times2$ quaternionic matrix $Q$; this is also briefly discussed in
\Ref{Corinne A. Manogue and Tevian Dray,
{\it Octonionic M\"obius Transformations},
Mod.\ Phys.\ Lett.\ {\bf A14}, 1243--1255 (1999).}%
.}

\Subsection{Decompositions}

One of the main results of \Eigen\ was that if $v$, $w$ are (normalized)
eigenvectors of the $2\times2$ octonionic Hermitian matrix $A$ corresponding
to different {\it real} eigenvalues $\lambda_v$, $\lambda_w$, then $A$ can be
expanded as
$$A = \lambda_v v v^\dagger + \lambda_w w w^\dagger \Eqno$$\label\Decomp
Furthermore, $v$ and $w$ are automatically orthogonal in the generalized sense
$$(v v^\dagger) \, w = 0 \Eqno$$\label\Ortho

We now ask whether a decomposition analogous to \Decomp\ exists when the
eigenvalues are not real.  Consider first the case of left eigenvalues, so
that $A\in\AAA$ and $v\in\VVV$.  If $v$ is given by \vform\ with
$v^\dagger v=1$, let
$$w:=\pmatrix{0& 1\cr 1& 0\cr} v = \pmatrix{y\cr x\cr} \in\VVV \Eqno$$
\label\Flip
This leads to $w^\dagger w=1$, and furthermore
$$v^\dagger w = \bar{x} \cdot \bar{y} = \bar{(y \cdot x)} = 0 \Eqno$$
by the definition of $\VVV$.  Provided
$$[i,x,y] = 0 \Eqno$$
so that the entire problem is quaternionic and hence associative, it turns out
that \Decomp\ does hold,
\Footnote{No parentheses are needed because of the assumed associativity.}
where $\lambda_v=\lambda$ is obtained by solving (either equation in) \Leq\
and $\lambda_w$ is obtained from $\lambda_v$ by interchanging $x$ and $y$.  We
illustrate this result by returning to the example \Unex, for which we obtain
the decomposition
$$\pmatrix{1&-i\cr \noalign{\smallskip} i&~~1} =
  {(1+j)\over2} \pmatrix{1\cr k} \pmatrix{1\cr k}^\dagger +
  {(1-j)\over2} \pmatrix{k\cr 1} \pmatrix{k\cr 1}^\dagger
  \Eqno$$
where the factor of $2$ is due to the normalization of the eigenvectors.

The above construction fails if $[i,x,y]\ne0$.  Remarkably, a similar
construction still works for {\it right\/} eigenvalues!  Direct computation
establishes:

\Theorem{For any $A\in\AAA$ and (normalized) $v\in\VVV$ such that
$Av=v\lambda_v$.  Then $A$ can be expanded as
$$A = \lambda_v \> (v v^\dagger) + \lambda_w \> (w w^\dagger)
  \Eqno$$\label\Surprise
where $w$ is defined by \Flip\ and satisfies $Aw=w\lambda_w$.}

\noindent
As before, we can assume without loss of generality that $A$ is given by 
\Aform\ and that $v$ and $w$ are given by \GenEx.  Returning to our example
\SpecEx\ yields the explicit decomposition
$$\pmatrix{~~1&i\cr \noalign{\smallskip} -i&1\cr} =
  {(1+k\ell)\over2}
	\left( \pmatrix{j\cr \ell} \pmatrix{j\cr \ell}^\dagger \right) +
  {(1-k\ell)\over2}
	\left( \pmatrix{\ell\cr j} \pmatrix{\ell\cr j}^\dagger \right)
  \Eqno$$

While it is true that
$$(v v^\dagger) \, w = v \, (v^\dagger w) \Eqno$$\label\VAssoc
for any $v$, $w$ related by \Flip\ (but not necessarily in $\VVV$), the
decomposition \Surprise\ is surprising because the eigenvalues $\lambda_v$,
$\lambda_w$ do {\it not\/} commute or associate with the remaining terms.
Specifically, although \VAssoc\ is zero here, we have
$$\Big(\lambda_v (v v^\dagger) \Big) \, w \ne 0 \Eqno$$
Remarkably, there is another decomposition theorem which does not have this
problem.  Direct computation establishes:

\Theorem{Given $A\in\AAA$ and (normalized) $v\in\VVV$ such that
$Av=v\lambda_v$.  Then $A$ can be expanded as
$$A =  (v \lambda_v) \, v^\dagger + (w \lambda_w) \, w^\dagger
  \Eqno$$\label\SurpriseII
where $w$ is defined by \Flip\ and satisfies $Aw=w\lambda_w$.}

\noindent
The decomposition \SurpriseII\ is less surprising than \Surprise\ when one
realises that orthogonality in the form
$$\Big( (v \lambda) \, v^\dagger \Big) \, w
  = (v \lambda) \, (v^\dagger w)
  = 0 \Eqno$$\label\NewOrtho
holds for {\it any\/} $\lambda\in\OO$ and $v,w\in\VVV$ satisfying \Flip.
Furthermore, since
$$\Big( (v \lambda) \, v^\dagger \Big) \, v
  = (v \lambda) \, (v^\dagger v) \Eqno$$\label\App
for {\it any\/} $v$ and $\lambda$ (see the Appendix), one can use the
decomposition \SurpriseII\ to construct $2\times2$ octonionic
matrices with arbitrary octonionic (right) eigenvalues.  However, such
matrices will {\it not\/} in general be Hermitian; this requires the imaginary
parts of the eigenvalues to be equal and opposite, as well as some further
restrictions on the eigenvalues (compare \GenEx).

\Section{SPIN}

In standard quantum theory, the infinitesimal generators $r_\alpha$ (with
$\alpha=x,y,z$) of angular momentum or spin are just the anti-Hermitian
matrices obtained by multiplying the Pauli matrices by the imaginary complex
unit (which for us is $\ell$, not $i$) and dividing by 2.  Explicitly, we have
$$r_x = {1\over2} \pmatrix{0& \ell\cr \noalign{\smallskip} \ell& 0\cr} \qquad
  r_y = {1\over2} \pmatrix{0& 1\cr \noalign{\smallskip} -1& 0\cr} \qquad
  r_z = {1\over2} \pmatrix{\ell& 0\cr \noalign{\smallskip} 0& -\ell\cr}
  \Eqno$$
where we have set $\hbar=1$.  One then normally multiplies by $-\ell$ to
obtain a description of the Lie algebra $su(2)$ in terms of Hermitian
matrices.

As discussed in [\DIM,\SPIN], however, in the octonionic setting care must be
taken with this last step, and we define instead the operators
$$L_\alpha \psi = - (r_\alpha \psi) \ell \Eqno$$
where $\psi$ is a 2-component {\it octonionic} column (representing a
Majorana-Weyl spinor in 10 spacetime dimensions).  The operators $L_\alpha$
are self-adjoint with respect to the inner product 
$$\langle \psi,\chi \rangle = \pi \!\left( \psi^\dagger\chi \right) \Eqno$$
where the map
$$\pi(q) = {1\over2} (q + \ell q \bar\ell) \Eqno$$
projects $\OO$ to a preferred complex subalgebra $\CC\subset\OO$.

Spin eigenstates are obtained as usual as the eigenvectors of $L_z$ with
eigenvalues $\pm{1\over2}$.  Particular attention is paid in [\DIM,\SPIN]
to the eigenstates
$$\eplus = \pmatrix{1\cr k\cr} \qquad
  \eminus = \pmatrix{-k\cr ~~1\cr}
  \Eqno$$
which were proposed as representing particles
\Footnote{These two papers used different conventions to distinguish particles
from antiparticles; we adopt the conventions used in \Spin.}
at rest with spin $\pm{1\over2}$, respectively.  Note that $\eplus$ and
$\eminus$ are orthogonal with respect to the above inner product, that is
$$\langle \eplus,\eminus \rangle = 0$$
We will therefore focus on these eigenstates, which have some extraordinary
properties.

Consider now the remaining spin operators $L_x$, $L_y$ acting on
these eigenstates.  We have
$$L_x \eplus = {1\over2} \pmatrix{-k\cr ~~1\cr}
  = \eplus \left( -{k\over2} \right) \Eqno$$
and
$$L_y \eplus = {1\over2} \pmatrix{-k\ell\cr \ell\cr}
  = \eplus \left( -{k\ell\over2} \right) \Eqno$$
with similar results holding for $\eminus$.  This illustrates the fact that
this {\it quaternionic} self-adjoint operator eigenvalue problem admits
eigenvalues which are not real.  More importantly, as claimed in [\DIM,\SPIN],
it shows that $\eplus$ is a {\it simultaneous\/} eigenvector of the 3
self-adjoint spin operators $L_x$, $L_y$, $L_z$!

This result could have significant implications for quantum mechanics.  In
this formulation, the inability to completely measure the spin state of a
particle, because the spin operators fail to commute, is thus ultimately due
to the fact that the {\it eigenvalues} don't commute.  Explicitly, we have
$$\eqalignno{
4L_x (L_y \eplus)
  &= 2L_x (-\eplus \,k\ell)
   = 2r_x (\eplus \,k\ell) \,\ell
   = -2r_x (\eplus \,\ell) \,k\ell
   = +(2L_x \eplus) \,k\ell \cr
  &= -\eplus \,k \, k\ell = + \eplus \,\ell \Eqalignno \cr
\noalign{\smallskip}
4L_y (L_x \eplus)
  &= -\eplus \,k\ell \, k = - \eplus \,\ell \Eqalignno
  }$$
which yields the usual commutation relation in the form
$$[L_x , L_y] \,\eplus
  = {1\over2} \, \eplus \,\ell = L_z \, \eplus \,\ell$$

Consider now the more general eigenstate
$$L_z \left( \eplus e^{\ell\theta} \right)
   = \left( \eplus e^{\ell\theta} \right) \left( {1\over2} \right) \Eqno$$
and note that
$$L_x \left( \eplus e^{\ell\theta} \right)
  = \left( \eplus e^{\ell\theta} \right)
    \left( -{k \, e^{2\ell\theta}\over2} \right) \qquad
L_y \left( \eplus e^{\ell\theta} \right)
  = \left( \eplus e^{\ell\theta} \right)
    \left( -{k\ell \, e^{2\ell\theta}\over2} \right) \qquad
\Eqno$$
so that the non-real eigenvalues depend on the phase.  It is intriguing to
speculate on whether the value of the non-real eigenvalues, which determine
the phase, can be used to specify (but not measure) the actual direction of
the spin, and whether this might shed some insight on basic properties of
quantum mechanics such as Bell's inequality.

\Section{\bold $3\times3$ OCTONIONIC HERMITIAN MATRICES}

The general $3\times3$ octonionic Hermitian matrix can be written
$$\AA = \pmatrix{p& a& \bar{b}\cr \noalign{\smallskip}
	\bar{a}& m& c\cr \noalign{\smallskip} b& \bar{c}& n\cr} \Eqno$$
\label\Three
with $p,m,n\in\RR$ and $a,b,c\in\OO$.  Remarkably, $\AA$ satisfies the
(natural generalization of the) usual characteristic equation for $3\times3$
matrices
\Ref{Hans Freudenthal,
{\it Zur Ebenen Oktavengeometrie},
Proc.\ Kon.\ Ned.\ Akad.\ Wet.\ {\bf A56}, 195--200 (1953).}%
.  Equally remarkably, the real eigenvalues do {\it not\/} satisfy this
equation, and, as shown originally by Ogievetsky
\Ref{
{\cyr O. V. Ogievetski\u\i},
{\cyr Kharakteristicheskoe U$\!$ravnenie dlya Matrits ${3\times3}$
		nad Oktavami}, 
Uspekhi Mat.\ Nauk {\bf 36}, 197--198 (1981); translated in:
O. V. Ogievetskii,
{\it The Characteristic Equation for ${3\times3}$ Matrices
	over Octaves},
Russian Math.\ Surveys {\bf 36}, 189--190 (1981).}%
, $\AA$ has 6, rather than 3, real eigenvalues.  As shown in \Eigen\ (see
also~\Okubo), the eigenvalues naturally belong to 2 distinct families, each
containing 3 real eigenvalues.  Furthermore, within each family, the
corresponding eigenvectors are orthogonal in the sense of \Ortho, and lead to
a natural decomposition of $\AA$ along the lines of \Decomp.


In our previous work \Eigen, we derived a generalized characteristic equation
satisfied by the real eigenvalues.  In this section, we summarize our
preliminary efforts to generalize this equation to non-real eigenvalues.  Since
the calculations are somewhat involved, and since we have so far been unable
to solve the resulting equations, we give here only the final results of those
calculations, the full details of which are available online together with
some examples and further discussion
\Ref{Tevian Dray, Jason Janesky, and Corinne A. Manogue,
{\it Some Properties of $3\times3$ Octonionic Hermitian Matrices with Non-Real
  Eigenvalues},
Oregon State University, 2000, 12 pages.
\hfill{\small({\tt http://xxx.lanl.gov/abs/math/0010255})}}\label\Online
.%

Inserting
$$v=\pmatrix{x\cr y\cr z\cr} \Eqno$$\label\vIII
into the (right) eigenvalue problem \Master\ leads, after considerable
manipulation \Online, to the generalized characteristic equation for the
eigenvalues $\lambda$ in the form
$$\eqalign{
  z \left(
	\lambda^3 - (\Tr \AA) \, \lambda^2 + \sigma(\AA) \, \lambda - \det \AA
	\right)
  &= b \Bigl( a (cz) \Bigr) + \bar{c} \Bigl( \bar{a} (\bar{b}z)	\Bigr) 
     - \Bigl( b(ac) + (\bar{c}\,\bar{a})\bar{b} \Bigr) z \cr
  &\qquad + b \, [a,y,\lambda] + [b,ay,\lambda]
          + [b,x,\lambda] \, (\lambda-m) \cr 
  &\qquad + \bar{c} \, [\bar{a},x,\lambda] + [\bar{c},\bar{a}x,\lambda]
          + [\bar{c},y,\lambda] \, (\lambda-p) \cr
  }\Eqno$$\label\CharIII
The term in parentheses on the left-hand-side is just the usual characteristic
equation; in particular, we have
$$\eqalign{
\Tr A &= p + m + n \cr
\sigma(A) &= pm + pn + mn - |a|^2 - |b|^2 - |c|^2 \cr
\det A &= pmn + b(ac) + \bar{b(ac)} - n|a|^2 - m|b|^2 - p|c|^2 \cr
  }\Eqno$$

If $\lambda$ is real, all the associators on the right-hand-side vanish, and
we recover the generalized characteristic equation given in \Eigen.  As shown
there, the requirement in that case that the right-hand-side be a real
multiple of $z$ (since the left-hand-side is) then constrains $z$, resulting
in precisely 2 values for that real multiple, and reducing \CharIII\ to 2
cubic equations, one for each family of real eigenvalues.

While we find the form of \CharIII\ attractive, as there are no extraneous
terms involving both $z$ and $\lambda$, we have so far been unable to further
simplify \CharIII\ when $\lambda$ is not real.  (It is straightforward to
eliminate one of $x$, $y$ from this equation, but not both.)

Another possible approach to finding the eigenvalues relies on the associator
identity
$$[v^\dagger,v,\lambda]
  := (v^\dagger v) \lambda - v^\dagger (v\lambda) \equiv 0 \Eqno$$\label\APP
which follows for {\it any} octonionic vector $v$ and $\lambda\in\OO$ by
alternativity, and which is further discussed in the Appendix.  If $v$ is a
normalized right eigenvector of $\AA$ with eigenvalue $\lambda$, then
$$v^\dagger (\AA v) = v^\dagger (v\lambda) = (v^\dagger v) \lambda = \lambda
  \Eqno$$\label\Jason
which yields an equation for $\lambda$ in terms of $\AA$ and the components of
$v$.  After some rearrangment, further details of which can again be found in
\Online, one obtains
$$\Re(\lambda) = {x \cdot (ay) + z \cdot (bx) + p |x|^2 \over |x|^2} \Eqno$$
(and similar expressions obtained by cyclic permutation) and
$$\Im(\lambda) = [x,a,y] + [z,b,x] + [y,c,z] \Eqno$$
We had hoped to use these various expressions to impose conditions on $\AA$
which would in turn enable us to solve for $\lambda$, but have not yet found a
way to do so.

\Section{DISCUSSION}

As pointed out in \Eigen, the orthonormality relation \Ortho\ is equivalent to
assuming that 
$$vv^\dagger + ... + ww^\dagger=I \Eqno$$
If we define a matrix $U$ whose columns are just $v,...,w$, then this
statement is equivalent to
$$U U^\dagger = I \Eqno$$\label\UOrtho
Furthermore, the eigenvalue equation \Master\ can now be rewritten in the form
$$A U = U D \Eqno$$ \label\Matrix
where $D$ is a diagonal matrix whose entries are the eigenvalues.
Decompositions of the form \SurpriseII\ now take the form
$$A = (UD) U^\dagger \Eqno$$
and multiplication of \Matrix\ on the right by $U^\dagger$ shows that
$$(AU) U^\dagger = (UD) U^\dagger = A = A (UU^\dagger) \Eqno$$
Thus, just as in \Eigen, decompositions of the form \SurpriseII\ can be viewed
as the assertion of associativity
$$(AU) U^\dagger = A (UU^\dagger) \Eqno$$\label\UAssoc
For non-real eigenvalues, we know of no way to express decompositions of the
form \Decomp\ in similar language, which leads us to suspect that \SurpriseII\
is more fundamental.  We further conjecture that the correct notion of
orthogonality is \NewOrtho, not \Ortho, to which it of course reduces if the
eigenvalues are real.  In any case, it is intriguing that this notion of
orthogonality can be written as
$$\Big( (\AA v) \, v^\dagger \Big) \, w = 0 \Eqno$$
which explicitly involves $\AA$.

Putting these ideas together, it would be natural to conjecture that {\it all}
eigenvectors of a $3\times3$ octonionic Hermitian matrix come in families of
3, which form a decomposition in the sense that \UAssoc\ is satisfied, and
which are orthogonal in the sense of \NewOrtho.  However, based on examples,
this conjecture appears to be false \Online.

There is, however, another intriguing possibility.  The examples considered in
\Online\ suggest that the eigenvectors of $3\times3$ octonionic Hermitian
matrices may come in sets of 6 (or more), rather than in sets of 3.  This
would fit nicely with our recent result with Okubo
\Ref{Tevian Dray, Corinne A. Manogue, and Susumu Okubo, in preparation.}
that, for real eigenvalues, it takes all 6 eigenvectors in order to decompose
an arbitrary vector into a linear combination of eigenvectors, despite the
fact that only 3 eigenvectors are needed to decompose the original matrix.

We therefore conjecture that, for any $3\times3$ octonionic Hermitian matrix,
\UAssoc\ should hold when suitably rewritten for a set consisting of $n$
eigenvectors, where $n$ presumably divides 24, the number of (real)
independent eigenvectors with real eigenvalues.  Whether or in what form
orthogonality would hold in such a context is an interesting open question.

\bigskip
\leftline{\bf APPENDIX}
\SUBSECNO=0

In deriving the foregoing results, we have made use of various associator
identities involving octonionic vectors, such as \App\ and \APP.  In this
appendix, we derive several such identities, including these two.  As we then
show, an important application of these identities is to give a particularly
elegant derivation of the so-called ``\Rule'' needed for supersymmetry.

\Subsection{Vector Associators}

Let $U$, $V$, $W$ be arbitrary octonionic vectors, i.e.\ $1\times n$
octonionic matrices.  Define the {\it vector associator} via
$$[U,V,W] := (U V^\dagger) W - U (V^\dagger W) \Eqno$$
Then direct computation using alternativity establishes
$$[W,V,V] \equiv 0 \Eqno$$\label\Master
which was used in \App.  Setting $W=V$ yields
$$[V,V,V] \equiv 0 \Eqno$$\label\Cube

Furthermore, essentially the same argument also establishes a similar formula
when $W$ is replaced by an octonionic scalar $\lambda$, so that
$$[\lambda,V,V] = 0 \Eqno$$
where we have implicitly defined yet another associator, namely
$$[\lambda,V,W] := (\lambda V^\dagger) W - \lambda (V^\dagger W) \Eqno$$
An interesting consequence of this result is the Hermitian conjugate relation
$$V^\dagger (V \bar\lambda) = (V^\dagger V) \bar\lambda \Eqno$$
which is equivalent to \APP.

We can polarize \Master\ to obtain
$$[U,V,W]+[U,W,V] \equiv 0 \Eqno$$\label\LemmaIII
a special case of which is
$$[V,V,W]+[V,W,V] \equiv 0 \Eqno$$
obtained by setting $U=V$.  A further special case of \LemmaIII\ is
$$[U,V,W]+[U,W,V]+[V,W,U]+[V,U,W]+[W,U,V]+[W,V,U] \equiv 0
  \Eqno$$\label\Six
obtained by adding cyclic permutations of \LemmaIII, or alternatively, without
requiring \Master, by repeated polarization of \Cube.

\Subsection{The 3-$\bold\Psi$'s Rule}

An essential ingredient in the construction of the Green-Schwarz superstring
[\MultiRef{M. B. Green, J. H. Schwarz,
{\it Covariant description of superstrings},
Phys.\ Lett.\ {\bf B136}, 367 (1984)},\Label\GS
\MultiRef{M. B. Green, J. H. Schwarz, and E. Witten,
{\bf Superstring Theory},
Cambridge University Press, Cambridge, 1987.}] \Label\GSB
is the spinor identity
$$\epsilon^{klm} \gamma^\mu \Psi_k \bar\Psi_l \gamma_\mu \Psi_m = 0 \Eqno$$
\label\Super
for anticommuting spinors $\Psi_k$, $\Psi_l$, $\Psi_m$, where $\epsilon^{klm}$
indicates total antisymmetrization.  This identity can be viewed as a special
case of a Fierz rearrangement.  An analogous identity holds for commuting
spinors $\Psi$, namely
$$\gamma^\mu \Psi \bar\Psi \gamma_\mu \Psi = 0 \Eqno$$\label\Commute
To our knowledge, Schray [\Thesis,\SCHRAY] was the first to formally refer to
\Super\ as the \Rule; we extend this usage to the commuting case \Commute.

It is well-known that the \Rule\ holds for Majorana spinors in 3 dimensions,
Majorana or Weyl spinors in 4 dimensions, Weyl spinors in 6 dimensions, and
Majorana-Weyl spinors in 10 dimensions.  Thus, the Green-Schwarz superstring
exists only in those cases [\GS,\GSB].  As was shown by Fairlie and Manogue
\Fairlie, the \Rule\ in all these cases is equivalent to an identity on the
$\gamma$-matrices, which holds automatically for the natural representation of
the $\gamma$-matrices in terms of the 4 division algebras $\RR,\CC,\HH,\OO$,
corresponding precisely to the above 4 types of spinors.  Manogue and Sudbery
\Sudbery\ then showed how to rewrite these spinor expressions in terms of
$2\times2$ matrices over the appropriate division algebra, thus eliminating
the $\gamma$-matrices completely.

In the commuting case, the (unpolarized) \Rule\ can be written in terms of a
2-component octonionic ``vector'' (really a spinor) $V$ as
\Schray\
$$\tilde{(VV^\dagger)} \, V = 0 \Eqno$$
where
$$\tilde{A} := A - \Tr{A} \Eqno$$
corresponding to time reversal.  It is straightforward to check that this
equation holds by alternativity.  This can also be seen by using the identity
$$\Tr(V V^\dagger) = \Tr(V^\dagger V) = V^\dagger V \in \RR \Eqno$$
to write
$$\tilde{(VV^\dagger)} \, V
  = (VV^\dagger - V^\dagger V) \, V
  = (VV^\dagger) V - V (V^\dagger V)
  = [V,V,V]
 \Eqno$$
and we see that the \Rule\ for commuting spinors is just \Cube.

Equivalently, we can rewrite the (polarized) \Rule\ in terms of
2-component ``vectors'' $U$, $V$, $W$ as
$$(\tilde{UV^\dagger}+\tilde{VU^\dagger}) \, W 
  + (\tilde{VW^\dagger}+\tilde{WV^\dagger}) \, U 
  + (\tilde{WU^\dagger}+\tilde{UW^\dagger}) \, V
  = 0 \Eqno$$\label\Three
Using the identity
$$\Tr{(U V^\dagger + V U^\dagger)}
  = \Tr{(V^\dagger U + U^\dagger V)}
  = V^\dagger U + U^\dagger V
  \Eqno$$
we have
$$\tilde{UV^\dagger}+\tilde{VU^\dagger}
  = (U V^\dagger + V U^\dagger) - (U^\dagger V + V^\dagger U) \Eqno$$
where the last term is Hermitian and hence real.  We thus have
$$(\tilde{UV^\dagger}+\tilde{VU^\dagger}) W
  = (U V^\dagger + V U^\dagger) W - W (U^\dagger V + V^\dagger U) \Eqno$$
Using this 3 times shows that the \Rule\ \Three\ is precisely the
same as \Six.

An analogous argument can be given for anticommuting spinors; this is
essentially the approach used in \Schray.  Combining these results, the \Rule\
can be written without $\gamma$-matrices in terms of 2-component octonionic
spinors $\psi_\alpha$ as
$$[\psi_1,\psi_2,\psi_3] \pm [\psi_1,\psi_3,\psi_2]
  + [\psi_2,\psi_3,\psi_1] \pm [\psi_2,\psi_1,\psi_3]
  + [\psi_3,\psi_1,\psi_2] \pm [\psi_3,\psi_2,\psi_1]
  \equiv 0 \Eqno$$
for both commuting ($+$) or anticommuting ($-$) spinors.  Both of these
identities follow from the identity \Cube\ applied to 2-component octonionic
vectors, which is a special case of the more general identity \Master\ which
holds for octonionic vectors of arbitrary rank.

\References

\bye